\documentstyle[12pt,fleqn,parskip]{article}
\newcommand{\tn}{\otimes}
\newcommand{\tnr}{\mbox{\scriptsize $\bigcirc$}\!\!\!\!\!\!\:
\raisebox{0pt}{\tiny $\rr$}}

\renewcommand{\t}{\mbox{\sc t}}
\newcommand{\r}{\mbox{\sc r}}
\newcommand{\id}{\mbox{\rm id}}
\newcommand{\rr}{{\cal R}}
\newcommand{\ff}{{\cal F}}
\renewcommand{\aa}{{\cal A}}
\newcommand{\uu}{{\cal U}}
\newcommand{\bu}{\mbox{\boldmath ${\cal U}$}}
\newcommand{\ba}{\mbox{\boldmath ${\cal A}$}}
\newcommand{\bm}{\mbox{\boldmath $m$}}

\newcommand{\bD}{\mbox{\boldmath $\Delta$}}
\newcommand{\bep}{\mbox{\boldmath $\varepsilon$}}
\newcommand{\bS}{\mbox{\boldmath $S$}}
\newcommand{\br}{\mbox{\boldmath ${\cal R}$}}
\newcommand{\bast}{{\mbox{\boldmath $\ast$}}}
\newcommand{\eins}{1\hspace{-3pt}{\rm I}}
\newcommand{\cc}{{\rm C}\!\!\!\mbox{\small\sf l}\;}
\newcommand{\beq}{\begin{equation}}
\newcommand{\eeq}{\end{equation}}
\newcommand{\bgar}{\begin{array}}
\newcommand{\ear}{\end{array}}
\newcommand{\con}{\begin{array}[b]{c}
\mbox{\footnotesize $\;\leftharpoonup$}\\}
\newcommand{\ccon}{\begin{array}[b]{c}
\mbox{\footnotesize $\;\rightharpoonup$}\\}

\newcommand{\crr}{\mbox{\asn
$\!\!\begin{array}[b]{c}\mbox{\footnotesize $\:\!\rightharpoonup$}\\
\rr\ear\!\!$\ase}}
\newcommand{\ccf}{\mbox{\asn
$\!\!\begin{array}[b]{c}\mbox{\footnotesize $\,\,\leftharpoonup$}\\
\ff\ear\!\!$\ase}}
\newcommand{\cf}{\mbox{\asn
$\!\!\begin{array}[b]{c}\mbox{\footnotesize $\,\,\rightharpoonup$}\\
\ff\ear\!\!$\ase}}
\newcommand{\ccv}{\mbox{\asn
$\!\!\begin{array}[b]{c}\mbox{\footnotesize $\:\!\leftharpoonup$}\\
v\ear\!\!$\ase}}
\newcommand{\cv}{\mbox{\asn
$\!\!\begin{array}[b]{c}\mbox{\footnotesize $\:\!\rightharpoonup$}\\
v\ear\!\!$\ase}}
\newcommand{\ccz}{\mbox{\asn
$\!\!\begin{array}[b]{c}\mbox{\footnotesize $\;\!\leftharpoonup$}\\
z\ear\!\!$\ase}}
\newcommand{\cz}{\mbox{\asn
$\!\!\begin{array}[b]{c}\mbox{\footnotesize $\,\rightharpoonup$}\\
z\ear\!\!$\ase}}

\begin{document}

\newcommand{\ase}{\renewcommand{\arraystretch}{1.4}}
\newcommand{\asn}{\renewcommand{\arraystretch}{0}}

\ase
\thispagestyle{empty}
\begin{flushright}
LMU--TPW 94--20 \\ q-alg/9412001 \\ November 1994
\end{flushright}
\vskip 1cm
\bigskip\bigskip
\begin{center}
{\bf\LARGE{Self-Diagonal Tensor Powers of}}
\end{center}

\begin{center}
{\bf\LARGE{Quantum Groups and R-Matrices for}}
\end{center}

\begin{center}
{\bf\LARGE{Tensor Products of Representations}}
\end{center}
\vskip 1.5truecm
\begin{center}{\bf Ralf A. Engeldinger}\footnote{supported by
Graduiertenkolleg {\sc Mathematik im Bereich ihrer Wechsel\-wirkung mit
der Physik} at Mathematisches Institut der Ludwig-Maximilians-Universit\"at
M\"unchen}
\\
\end{center}
\begin{center}
{\bf }
\vskip10mm
Sektion Physik der Ludwig-Maximilians-Universit\"at M\"unchen\\
Lehrstuhl Professor Wess\\ Theresienstra\ss e 37, D-80333 M\"unchen,
Federal Republic of Germany\\

\vskip5mm
e-mail: engeldin@lswes8.ls-wess.physik.uni-muenchen.de
\end{center}
\vskip 1cm
\begin{abstract}
Twisted tensor powers of quasitriangular Hopf algebras with
diagonal sub-Hopf-algebras (self-diagonal tensor powers) are introduced
together with their duals and their mutual $\ast$-structures as
generalizations of the Drinfel'd double as given by Reshetikhin and
Semenov-Tian-Shansky. R-Matrices for tensor products of representations are
derived.
\end{abstract}

\newpage
\section{Introduction}
In [CEJSZ] the Drinfel'd double of a quasitriangular Hopf algebra (QTHA)
$(\uu,\rr)$ as
a twisted tensor product of two copies of $\uu$ (originally given by
Reshetikhin and Semenov-Tian-Shansky [ReSe])
was recovered as the result of an
attempt to define a tensor product of two copies of $\uu$ in such a way
that the image of the coproduct of $\uu$ is a
(so-called diagonal) sub-Hopf algebra (which for the
standard tensor product Hopf algebra is not the case).

Tensor products of representations of Hopf algebras are defined via the
coproduct. This way one immediately obtains R-matrices for tensor products of
representations of quasitriangular Hopf algebras from the universal R-matrix.
But the R-matrices obtained this way are not the
only possible ones.
Lorek, Schmidke and Wess [LSW] constructed all R-matrices for the
$[3]\oplus [1]$--representation of $U_qsu(2)$: the standard
$SO_{q^2}(3)$-R-matrix
and the two $SO_q(1,3)$-R-matrices. The former is the one obtainable
by applying the coproduct on the universal R-matrix of $U_qsu(2)$. The latter
ones, on the other hand, were shown in [CEJSZ] to naturally arise from the
universal R-matrices of the twisted tensor product of two copies of
$U_qsu(2)$.
It was this observation that inspired the present work.

Explicitly it looks
as follows. Using the product form of the representation
\[ [2]\times[2]=[3]\oplus [1] \]
and the notation
\[ \langle [i]\!\times\![j],h\rangle=\langle [i]\tn[j],\Delta(h)\rangle\]
we can write
{\samepage
\beq\label{so3} \langle ([2]\!\times\![2])\tn([2]\!\times\![2]),
\rr\rangle=R_{SO(3)}\oplus
1_{[3]\tn[1]}\oplus 1_{[1]\tn[3]}\oplus 1_{[1]\tn[1]} \eeq
\beq\label{lor1} \langle ([2]\!\tn\![2])\tn([2]\!\tn\![2]),
\rr^{-1}_{41}\rr_{13}\rr^{-1}_{42}\rr_{23}
\rangle = R_{I_{SO(1,3)}}\eeq
\beq\label{lor2} \langle ([2]\!\tn\![2])\tn([2]\!\tn\![2]),\rr^{-1}_{41}
\rr^{-1}_{31}
\rr^{-1}_{42}\rr_{23}\rangle = R_{II_{SO(1,3)}}\eeq}

In this paper the same
is done for the tensor product of an arbitrary number $s$ of copies of $\uu$.
Like in [CEJSZ] for the case $s=2$ the corresponding Hopf dual (for $s=2$ also
discussed in [Pod]) and the
respective $\ast$-structures are given. Hence self-diagonal tensor power
Hopf
algebras can be regarded as generalizations of the Drinfel'd double or
of complex quantum groups [DSWZ].
Following Wess' idea we then
derive R-matrices for tensor products of representations from the universal
R-matrices of quasitriangular self-diagonal tensor power Hopf algebras.

\section{Self-diagonal Tensor Power Hopf Algebras}
The tensor product of two copies of a quasitriangular Hopf algebra $(\uu,\rr)$
is again a QTHA in a natural way. The same is true for the tensor
product of an arbitrary number of such copies.

To number tensor factors of $\uu^{\otimes s}\otimes\uu^{\otimes s}$
we use both a natural notation
\[ (1_1,1_2,\ldots,1_s,2_1,\ldots,2_s) \]
and a `flattened' one
\[ (1,2,\ldots,s,s+1,\ldots,2s), \]
and additionally
\[ \bar 1=(1,\ldots,s)=(1_1,\ldots,1_s),\qquad \bar 2=(s+1,\ldots,2s)=
(2_1,\ldots,2_s). \]
\[ \t(a\otimes b)=b\otimes a\]
\renewcommand{\arraystretch}{.5}
\[ \t_i=\id^{\otimes i-1}\otimes\t\otimes\id^{\otimes 2s-i-1}\]
\beq \t^{\left(s\right)}_\Delta=
\begin{array}{c}\mbox{\scriptsize{$s\!-\!2$}}\\ \bigcirc\\
\mbox{\scriptsize $i\!=\!0$}
\end{array}\left(
\begin{array}{c}\mbox{\scriptsize $i$}\\ \bigcirc\\ \mbox{\scriptsize
$j\!=\!0$}
\end{array}\t_{s-i+2j}\right),\qquad
\t^{\left(s\right)}_m=\t^{\left(s\right)-1}_\Delta
=\begin{array}{c}\mbox{\scriptsize{$0$}}\\ \bigcirc\\
\mbox{\scriptsize $i\!=\!s\!-\!2$}
\end{array}\left(
\begin{array}{c}\mbox{\scriptsize $i$}\\ \bigcirc\\ \mbox{\scriptsize
$j\!=\!0$}
\end{array}\t_{s-i+2j}\right).\eeq
$\bigcirc$ denotes iterated concatenation of mappings ($\circ$).

\renewcommand{\arraystretch}{1.4}
Then the standard tensor power QTHA has the following form:
\beq\label{std}\bgar{l}\bar\uu\equiv\uu^{\tn s},\quad\bar m=m^{\tn s}\!\circ\!
\t^{\left(s\right)}_m,\quad
\bar\Delta=\t^{\left(s\right)}_\Delta\!\circ\!\Delta^{\tn s},\quad
\bar\varepsilon=\varepsilon^{\tn s},\quad\bar S=S^{\tn s},\\
\bar\rr_{\bar 1\bar 2}=\prod\limits_{i=1}^s\rr^{\left(i\right)}_{i,s+i}
=\prod\limits_{i=1}^s\rr^{\left(i\right)}_{1_i2_i} \ear\eeq
where for any $i:\quad\rr^{\left(i\right)}_{12}\in\{\rr_{12},\rr^{-1}_{21}\}$.

In straightforward generalization of the case $s\!=\!2$ we now introduce a
different QTHA structure on the same algebra which has a diagonal sub-Hopf
algebra isomorphic to $\uu$. It is obtained from the standard tensor power
QTHA through a twist [Dri]. We denote this self-diagonal tensor product by
{\scriptsize $\bigcirc$}$\!\!\!\!\!\:$\raisebox{0pt}{\tiny $\rr$}.
Introducing the following left action on a Hopf algebra:
\beq \cz(x)=z\,x\,z^{-1}\eeq
and
\renewcommand{\arraystretch}{.5}
\beq \r^{\ff}_k=\t_k\!\circ\!\crr^{\ff}_{k,k+1},\qquad
\r^{\left(s\right)}=
\begin{array}{c}\mbox{\scriptsize $s\!-\!1$}\\ \bigcirc\\
\mbox{\scriptsize $i\!=\!1$}
\end{array}\left(
\begin{array}{c}\mbox{\scriptsize $i$}\\ \bigcirc\\
\mbox{\scriptsize $j\!=\!1$}
\end{array}\r^{\ff}_{s-i-1+2j}\right)\eeq

\renewcommand{\arraystretch}{1.4}

\beq\label{F}\bgar{l} \ff=\prod\limits_{i=1}^{s-1}\left(\prod\limits_{j=1}^i
\rr^{\ff}_{s+j,s-i+j}\right)
=\prod\limits_{i=s-1}^1\left(\prod\limits_{j=1}^{s-i}\rr^{\ff}_{2_j1_{i+j}}
\right)\\ \\
v=\ff^{(\bar 1)}\bar S(\ff^{(\bar 2)})
=\left[\prod\limits_{i=1}^{s-1}\left(\prod\limits_{j=i+1}^s\rr^\ff_{ij}
\right)\right]^{-1}
,\quad
v^{-1}=\bar S^2(\ff^{(\bar 1)})\ff^{(\bar 2)}
\ear\eeq
where $\rr^{\ff}_{12}\in \{\rr_{12},\rr^{-1}_{21}\}$,
we define the self-diagonal tensor power QTHA by
\beq\bgar{l}\bu\equiv\uu^{\tnr s}=\uu^{\tn s}\mbox{ as vector space over }\cc
,\\
\bm=
\bar m,\quad\bD=\r^{(s)}\!\circ\!\Delta^{\tn s}=
\cf\!\circ\!\bar\Delta,\quad
\bep=\bar\varepsilon,\quad\bS=\cv\!\circ\!\bar S,\\ \br_{\bar 1\bar 2}=
\ff_{\bar 2\bar 1}\bar\rr_{\bar 1\bar 2}\ff^{-1}_{\bar 1\bar 2}.\ear\eeq

Since the definition of the standard tensor power Hopf algebra is
self-dual we can obtain the dual Hopf algebra (Hopf dual) of the tensor power
QTHA $(\uu^{\tn s},\bar\rr)$ of a QTHA $(\uu,\rr)$ as the tensor power
Hopf algebra $\aa^{\tn s}$ of the latter one's dual $\aa$. By dualization of
the twist transformation we get the dual of the self-diagonal tensor power QTHA
$(\uu^{\tnr s},\br)$ which we denote as $\aa^{\tnr s}$. (The symbol
{\scriptsize $\bigcirc$}$\!\!\!\!\!\:$\raisebox{0pt}{\tiny $\rr$} has
a different meaning here.) Introducing the following right action
on the dual Hopf algebra
\beq\ccz(a)=\langle a_{(\bar 1)},z\rangle\;a_{(\bar 2)}\;
\langle a_{(\bar 3)},z^{-1}\rangle\eeq
we define it to be
\beq\bgar{l}\ba\equiv\aa^{\tnr s}=\bar\aa \mbox{ as vector space over
$\cc$},\\
\bm=\bar m\circ\!\ccf
,\qquad\bD=\bar\Delta,\qquad
\bep=\bar\varepsilon,\qquad\bS=\bar S\circ\!\ccv.
\ear\eeq
A tensor power Hopf algebra allows various inequivalent $*$-structures (i.e.
not related via (co-)conjugation with an invertible element: $h^\dagger\neq
\cz(h^\ast)$
or $a^\dagger\neq
\ccz(a^\ast)$). The naive choice
is simply $\bar\ast=\ast^{\tn s}$. From this we can obtain
an inequivalent one for every permutation $\pi\in {\cal S}_s$ with
$\pi^2=\id$ :
$\bar\ast_\pi=\pi\!\circ\!\bar\ast$.

We now focus on the case
\beq \rr^{\ast\tn\ast}_{12}=\rr_{21}. \eeq

Denoting by $\iota$ the inversion in ${\cal S}_s$ this implies:
\beq\label{v} \iota(v^{\bar\ast})=v.\eeq

Let $u\in\aa\tn\cc^{n\times n}$ be a fundamental representation of $\uu$
(generating all irreducible ones) such that
$R_{12}=\langle u_1\tn u_2,\rr\rangle$ fulfills
\beq \overline{R^\top_{12}}=R_{21}. \eeq
In $\aa^{\tnr s}$ we have\footnote{$a\bullet b=\bm(a\tn b)$}
\renewcommand{\arraystretch}{2}
\beq\bgar{l} \langle b^\bast{}_{(\bar 1)}\tn
a^\bast{}_{(\bar 1)},\ff\rangle\;b^\bast{}_{(\bar 2)}a^\bast{}_{(\bar 2)}\;
\langle
b^\bast{}_{(\bar 3)}\tn a^\bast{}_{(\bar 3)},\ff^{-1}\rangle =b^\bast\bullet
a^\bast\\ \qquad=(a\bullet b)^\bast
=\overline{\langle a_{(\bar 1)}\tn b_{(\bar 1)},\ff\rangle}\;
(a_{(\bar 2)}b_{(\bar 2)})^\bast\;
\overline{\langle a_{(\bar 3)}\tn b_{(\bar 3)},\ff^{-1}\rangle}.\ear\eeq
\ase
With
\beq F=\prod_{i=1}^{s-1}\left(\prod_{j=1}^i R^\ff_{s+j,s-i+j}\right)
=\prod\limits_{i=s-1}^1\left(\prod\limits_{j=1}^{s-i}R^	\ff_{2_j1_{i+j}}
\right)\\
\eeq
this implies
$\overline{F_{\bar 1\bar 2}}=(\iota\tn\iota)(F_{\bar 1\bar 2})$ if we assume
$\bast=\bar\ast$. In general this is not true. However there is a different
choice which always works in our specified case:
\beq \bast=\iota\!\circ\!\bar\ast\qquad ( \mbox{on }\ba). \eeq
Because of eq.(\ref{v}) we can define (cf.[CEJSZ])
\beq \bast=\cv\!\circ\!\iota\!\circ\!\bar\ast\qquad ( \mbox{on }\bu). \eeq
The natural generators of $\aa^{\tnr s}$ are
\beq u^{(i)}=\eins^{\tn i-1}\tn u\tn\eins^{\tn s-i}\eeq
while those of $\uu^{\tnr s}$ are as usual the corresponding
semirepresentations of $\br$ \linebreak and $\br^{-1}$
\beq L^{+(i)}=\langle\,\cdot\tn u^{(i)},\br\rangle \eeq
\beq L^{-(i)}=\langle u^{(i)}\tn\cdot\,,\br^{-1}\rangle. \eeq
They consist of $s+1$ different matrices of
generators\footnote{$(\ell^+_1\tn\ell^+_1)^m_n=\ell^+{}_k^m\tn\ell^+{}_n^k$}
(cf.[CEJSZ]):
\beq \ell^{(j)}_1=\bigotimes\limits_{k=1}^{s-j}\ell^+_1\tn
\bigotimes\limits_{k=1}^j\ell^-_1,\qquad j=0,\ldots,s\qquad\mbox{if}\quad
\rr^\ff=\rr\eeq
or
\beq \ell^{(j)}_1=\bigotimes\limits_{k=1}^{s-j}\ell^-_1\tn
\bigotimes\limits_{k=1}^j\ell^+_1,\qquad j=0,\ldots,s\qquad\mbox{if}\quad
\rr^\ff=\rr^{-1}_{21}\eeq
where$\quad\ell^+=\langle\,\cdot\tn u,\rr\rangle\quad$and$\quad
\ell^-=\langle u\tn\cdot\,,\rr^{-1}\rangle$.

We close this section with the following observation:

$\Delta^{s-1}$ is a $\ast$-algebra isomorphism\footnote{$\Delta^1=\Delta,
\enspace\Delta^n=(\Delta^{n-1}\tn\id)\!\circ\!\Delta$}
 and as such canonically gives rise
to an isomorphism of quasitriangular $\ast$-Hopf algebras. Its image is a
quasitriangular sub-$\ast$-Hopf algebra of $\uu^{\tnr s}$. Therefore
$(\uu^{\tnr s},\br)$ is called {\sl quasitriangular self-diagonal tensor power
$\ast$-Hopf algebra}.

\section{R-Matrices for Tensor Products of \hspace{10cm}$\;$ Representations}
Tensor products of representations naturally live on the corresponding tensor
power of the algebra. This way a coproduct on the algebra endows the category
of representations of the algebra with a monoidal structure. With the image
of $\Delta^{s-1}$ being a sub-$\ast$-Hopf algebra of $\uu^{\tnr s}$ we can
take advantage of the fact that the latter one is equipped with
quasitriangular structures that are not inherited directly from $(\uu,\rr)$
via the isomorphism. They give rise to R-matrices for tensor products of
representations of $(\uu,\rr)$ that cannot be obtained by simply applying
$\Delta^{s-1}\tn\Delta^{s-1}$ on $\rr$ or $\rr^{-1}_{21}$.
To the latter ones
we will refer as the {\sl contracted}
tensor product (cf. eq. \ref{so3})
and to the former ones as the
{\sl uncontracted} one (cf. eqs. \ref{lor1},\ref{lor2}).
Note that reordering factors within a contracted
tensor product leads to an equivalent R-matrix because of quasitriangularity.

Within tensor products of arbitrary many representations every single tensor
product can be realized as either a contracted or an uncontracted one.
After choosing a particular contraction of a tensor product
of representations we determine the remaining number of tensor factors $s$.
Now we select one of the $2^s$ possible universal R-matrices $\br$ of
$\uu^{\tnr s}$. Pairing it with two copies of the selected contraction gives us
an R-matrix for our original (precontracted) representation. (Of course this
procedure also works if we take two different representations (with the same
$s$) for the two components of $\br$.)

In the case of the tensor product of three representations we have the
following possibilities (neglecting permutations):
\beq\bgar{ll}s=1 & \langle([\enspace]\!\times\![\enspace]\!\times\![\enspace])
\tn([\enspace]\!\times\![\enspace]
\!\times\![\enspace]),\br\rangle,\qquad\qquad\mbox{here: }\br=\rr^{(1)}\\
s=2 &
\langle([\enspace]\!\tn\!([\enspace]\!\times\![\enspace]))\tn([\enspace]\!
\tn\!([\enspace]\!\times\![\enspace])),\br\rangle \\
& \langle(([\enspace]\!\times\![\enspace])\!\tn\![\enspace])\tn(([\enspace]\!
\times\![\enspace])\!\tn\!
[\enspace]),\br\rangle \\
s=3 & \langle([\enspace]\!\tn\![\enspace]\!\tn\![\enspace])\tn
([\enspace]\!\tn\![\enspace]\!\tn\![\enspace]),\br\rangle \\ \ear\eeq

Being interested in differential calculi on quantum spaces one has to
determine the eigenvalues of $\hat R=PR$. As shown by Wess and Zumino
[Wess, WeZu]
such a differential calculus can be defined if $\hat R$ has only one negative
eigenvalue (belonging to the antisymmetric projector). The R-matrices arising
from our procedure don't have this quality for $s>2$ as is easily seen by a
combinatorial argument. In the case $s=3$, e.g., we have four contributions
($++-\; ,\;+-+\; ,\;-++\; ,\;---$). The
eigenvalues depend only on the $\rr^{(i)}$ (in eq.(\ref{std}), and not on
$\rr^{\ff}$ in eq.(\ref{F})).
To keep the antisymmetric projector from splitting up it is necessary that
all three corresponding matrices $\hat R^{(i)}$ have exactly two eigenvalues
with $\lambda^{(i)}_-=-\lambda^{(i)}_+$. In all cases of interest this is not
fulfilled.

{\sl The first part of this work is based on [CEJSZ]. I appreciate the
attention,
support and  critique from K. F\" orger, P. Schupp,  S. Theisen and J. Wess.
I'm especially grateful to J. Wess for directing my attention to R-matrices
of the kind discussed here. I owe ongoing inspiration to works by V.G.
Drinfel'd.}

\end{document}